\documentclass [11pt,oneside]{amsart}
\usepackage{amsmath}

\usepackage {amsfonts}
\usepackage {amsmath}
\usepackage {amsthm}
\usepackage {amssymb}
\usepackage {framed}
\usepackage {amsxtra}
\usepackage {enumerate}
\usepackage {graphicx}
\usepackage{color}
\usepackage{graphicx}
\usepackage{graphicx}
\usepackage{wrapfig}
\usepackage{lscape}
\usepackage{rotating}
\usepackage{epstopdf}
\usepackage{pdflscape}
\usepackage{setspace}
\usepackage[left=0.86in, right=.86 in, top=1.2 in, footskip=0.5 in, bottom=0.9 in]{geometry}
\usepackage{adjustbox}
\usepackage[skip=0pt]{caption} 
\usepackage{tikz}
\usepackage{hyperref}
\usepackage{amsaddr}

\makeatletter
\renewcommand{\email}[2][]{%
  \ifx\emails\@empty\relax\else{\g@addto@macro\emails{,\space}}\fi%
  \@ifnotempty{#1}{\g@addto@macro\emails{\textrm{(#1)}\space}}%
  \g@addto@macro\emails{#2}%
}
\makeatother

\makeatletter

\theoremstyle{definition}
\newtheorem{df}{Definition} [section]

\theoremstyle{plain}
\newtheorem{thm}[df]{Theorem}

\newtheorem{lemma}[df]{Lemma}

\newtheorem{problem}[df]{Problem}

\newtheorem{conj}[df]{Conjecture}
\newtheorem{opr}[df]{Open Problem}

\title{The Hadwiger-Nelson problem with two forbidden distances}
\author{Geoffrey Exoo}
\address{Department of Mathematics and Computer Science, Indiana State University, Terre Haute, IN 47809}
\email{ge@cs.indstate.edu}

\author{Dan Ismailescu}
\address{Mathematics Department, Hofstra University, Hempstead, NY 11549}
\email{dan.p.ismailescu@hofstra.edu}

\begin{document}

\begin{abstract}
In 1950 Edward Nelson asked the following simple-sounding question:

\emph{How many colors are needed to color the Euclidean plane $\mathbb{E}^2$ such that no two points distance $1$ apart are identically colored?}

We say that $1$ is a \emph{forbidden} distance. For many years, we only knew that the answer was $4$, $5$, $6$, or $7$.
In a recent breakthrough, de Grey \cite{degrey} proved that at least five colors are necessary.

In this paper we consider a related problem in which we require \emph{two} forbidden distances, $1$ and $d$.
In other words, for a given positive number $d\neq 1$, how many colors are needed to color the plane such that no
two points distance $1$ \underline{or} $d$ apart are assigned the same color?
We find several values of $d$, for which the answer to the previous question is
at least $5$. These results and graphs may be useful in constructing simpler $5$-chromatic unit distance graphs.

\end{abstract}

\maketitle
\pagenumbering{arabic}
\newpage
\section{\bf Introduction}

Let $(\mathbb{E}^2,\|\cdot\|)$ be the $2$-dimensional space with the usual
Euclidean norm, and let $\mathcal{D}$ be a subset of $(0,\infty)$.
The \emph{distance graph} defined by $\mathcal{D}$ is the graph
$G(\mathbb{E}^2,\mathcal{D})$ whose vertices are points in $\mathbb{E}^2$ and
whose edges are the pairs of points $\{x,y\}$ such that $\|x-y\|\in \mathcal{D}$.

The graph $G$ has \emph{chromatic number} $\chi(\mathbb{E}^2,\mathcal{D})=k \in \mathbf{N}$ if there exists a $k$-coloring of $G$, i.e. a function from $\mathbb{E}^2$ to $\{1,2,\ldots, k\}$ that maps adjacent vertices to different values, and $k$ is minimal with this property.

A long-standing open problem is to determine the chromatic number of $\mathbb{E}^2$ when the distance set $\mathcal{D}$ reduces to only one element, say $\mathcal{D}=\{1\}$.
\begin{opr}{\bf The Hadwiger-Nelson Question.}
How many colors are needed to color the points of the plane such that no two points distance $1$ apart are assigned the same color?
\end{opr}
According to Soifer \cite{soifer}, the question was raised by Edward Nelson in 1950. It is easy to show that
\begin{equation}\label{bounds}
4\le \chi(\mathbb{E}^2,\{1\})\le 7.
\end{equation}
The upper bound follows from the existence of a tiling of the plane by regular hexagons, of diameter strictly between $2/\sqrt{7}$ and  $1$, which can be assigned seven colors in a periodic manner to form a $7$-coloring of the plane with no two identically colored points distance $1$ apart. This construction first appeared in a paper by Hadwiger \cite{hadwiger}.

\begin{figure}[htp]
\centering
\includegraphics[width=0.9\linewidth]{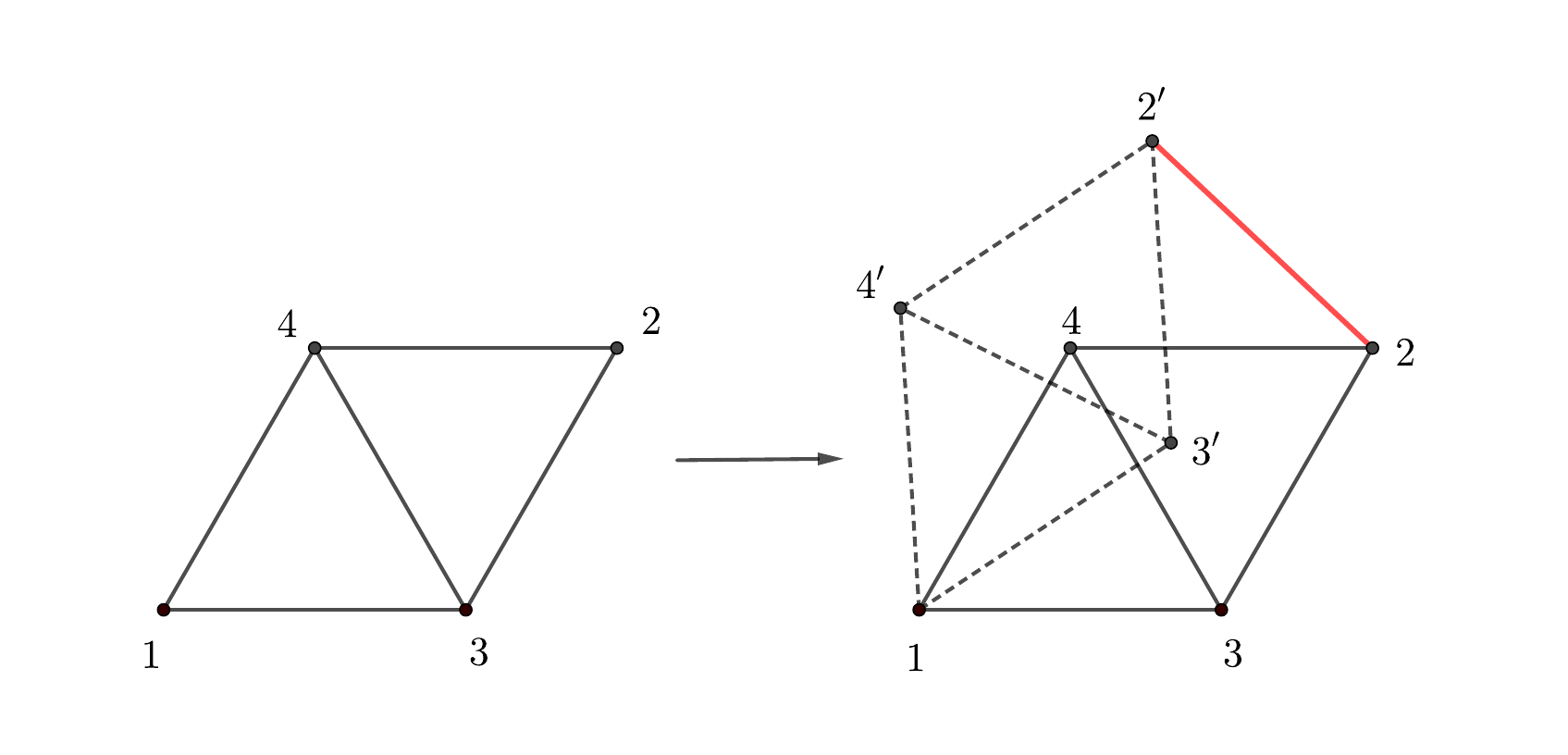}
\caption{Mosers' spindle construction}
\label{spindlefig}
\end{figure}

The lower bound is due to Leo and William Moser \cite{mosermoser}. Rotate a unit rhombus around one of the vertices of degree $2$ until the opposite vertex is at distance $1$ from its original position. The resulting $7$-vertex graph, called the \emph{Mosers's spindle}, has chromatic number $4$ - see figure \ref{spindlefig}.

Very recently, de Grey \cite{degrey} managed to improve the lower bound to $\chi(\mathbb{E}^2)\ge 5$. In an attempt to better understand the structure and properties
of $5$-chromatic unit distance graphs, it is of interest to consider that case when the set of distances $\mathcal{D}$ contains more than one element.

Exoo \cite{exoo} studied the case when the set is an entire interval, $\mathcal{D}=[1,d]$, and proved that
\begin{equation*}
\chi(\mathbb{E}^2,[1, 1.00853\ldots])\ge 5, \quad \text{and}\quad \chi(\mathbb{E}^2,[1, d])=7, \quad \text{for}\,\, 1.13475\ldots5 <d <1.38998\ldots.
\end{equation*}

The case when the distance set consists of all odd positive integers was
investigated by Ardal, Ma\v{n}uch, Rosenfeld, Shelah, and Stacho \cite{ardal}.
They proved that
\begin{equation*}
\chi(\mathbb{E}^2,\{1,3,5,7,9,11,\ldots\})\ge 5.
\end{equation*}

An interesting case occurs when $|\mathcal{D}|=2$, that is, when there are exactly two forbidden distances. We will be studying the following:

\begin{problem}\label{mainproblem}
Find values $d\neq 1$ such that $\chi(\mathbb{E}^2, \{1,d\})\ge 5$.
\end{problem}

Despite being a natural question, it appears that until very recently \cite{webpage} the problem had not been considered in this explicit form.
The only relevant result we could find is due to Katz, Krebs and Shaheen \cite{KKS}, who proved
\begin{equation*}
\chi(\mathbb{E}^2, \{1,\sqrt{2}\})\ge 5.
\end{equation*}
However, their initial motivation was different, as we will explain later on, when we will also provide a new proof of this result.

There are two good reasons we believe that Problem \ref{mainproblem} is worth looking into.
On one hand, it should be easier to prove that $\chi(\mathbb{E}^2, \{1,d\})\ge
5$ than to prove $\chi(\mathbb{E}^2, \{1\})\ge 5$.
On the other hand, as recently shown in \cite{EI}, one can prove results of the following type
\begin{thm}
There exist values of $d\neq 1$ such that
\begin{equation}
\text{If}\quad \chi(\mathbb{E}^2, \{1\})=4\quad \text{then}\quad \chi(\mathbb{E}^2, \{1,d\})=4.
\end{equation}
\end{thm}
The authors of \cite{EI} used $d=\sqrt{11/3}$ but most likely such an implication can be proved for other values of $d$.
They then proved that $\chi(\mathbb{E}^2, \{1,\sqrt{11/3}\})\ge 5$, by constructing a $5$-chromatic $\{1,\sqrt{11/3}\}$-graph, and therefore obtained that
$\chi(\mathbb{E}^2, \{1\})\ge 5$.

Before we proceed to our results, we state the following generalization of Mosers' spindle idea.
\begin{thm}{\bf The Spindle Method.}
Let $G$ be a finite graph with vertex set $V=\{1,2,\ldots, n\}$ and edge set $E$. Assume that the chromatic number of $G$, $\chi(G)=k$, and that in every $k$-coloring of $G$,
vertices $1$ and $2$ are colored identically.

Let $G'$ be a copy of $G$ such that $1=1'$ and $2\neq 2'$. Then the chromatic number of the graph $H$ whose edge set is $E\cup E' \cup \{\{2,2'\}\}$ is $\ge k+1$.
\end{thm}
\begin{proof} Assume that $H$ is $k$-colorable and let $c$ be such a $k$-coloring. Then $c(1)=c(2)$ by our assumption about $G$. At the same time, since $G'$ is a copy of $G$, it follows that $c(1')=c(1)=c(2')$. But then vertices $2$ and $2'$ have the same color, which violates the condition that the endpoints of edge $\{2,2'\}$ must be colored differently.
\end{proof}

We will also use the following definition very frequently.
\begin{df}
For a given positive real number $d\neq 1$, a \emph{$\{1,d\}$-graph} is a finite graph whose vertices are points in the Euclidean plane $\mathbb{E}^2$, and whose edges are obtained by connecting two points whenever the distance between them is either $1$ or $d$.
\end{df}
It is clear that given $d\neq 1$, if one can find a $\{1, d\}$-graph whose chromatic number is $5$, then it would immediately follow that $\chi(\mathbb{E}^2,\{1,d\})\ge 5$.
Most of our paper is dedicated to constructing such graphs, for various values of $d$.

Also note that any $\{1,d\}$-graph can be transformed into a $\{1, 1/d\}$-graph by an appropriate scaling. For this particular reason, we will always consider only the case $d>1$.

Hence, we want to construct $\{1,d\}$-graphs that have chromatic number $5$. Ideally, these graphs would have not too many vertices, as computing the chromatic number of large graphs
is computationally expensive.

The natural idea is to require such a graph to contain many small subgraphs that are $4$-chromatic. The smallest $4$-chromatic graph is of course $K_4$, the complete graph on $4$ vertices. The question now becomes, for what values of $d\neq 1$ can $K_4$ be represented as a $\{1,d\}$-graph? Luckily for us, this problem was already solved by Erd\H{o}s and Kelly \cite{EK}, and later by Einhorn and Schoenberg \cite{ES}.

\begin{thm}\label{twodistance}\cite{EK, ES}
The only values $d>1$ for which the complete graph $K_4$ can be embedded as a $\{1,d\}$-graph are: $d=(\sqrt{5}+1)/2$, $d=\sqrt{3}$, $d=(\sqrt{6}+\sqrt{2})/2$, and $d=\sqrt{2}$.
\end{thm}
\begin{figure}[htp]
\centering
\includegraphics[width=\linewidth]{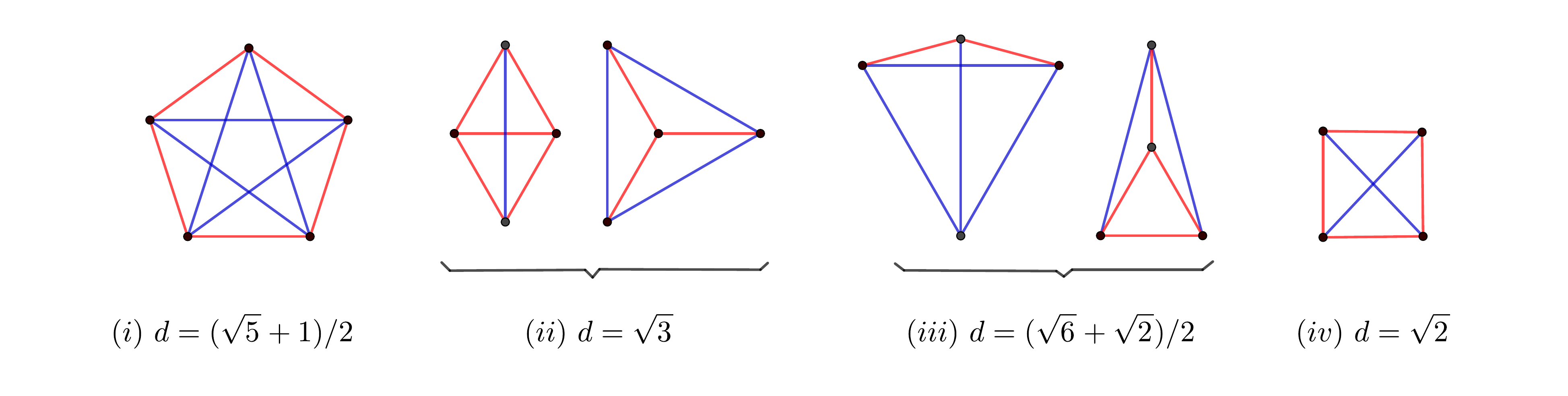}
\caption{Complete graphs realized as $\{1,d\}$-graphs}
\label{k4k5fig}
\end{figure}
Here and later we make the convention of representing the edges of length $1$ in red and the edges of length $d$ in blue - see figure \ref{k4k5fig}.
Note that if $d=(\sqrt{5}+1)/2$, then $K_5$, (the complete graph on $5$ vertices) can be expressed as a $\{1, d\}$-graph. Since $\chi(K_5)=5$, we immediately have that
\begin{thm}\label{root5thm}
\begin{equation*}
\chi\left(\mathbb{E}^2, \left\{1, (\sqrt{5}+1)/2\right\}\right)=\chi\left(\mathbb{E}^2, \left\{1, (\sqrt{5}-1)/2 \right\}\right)\ge 5.
\end{equation*}
\end{thm}
\newpage
\section{\bf Two simple cases: $d=\sqrt{3}$ and $d=(\sqrt{6}+\sqrt{2})/2$.}

If $d\in \{\sqrt{3}, (\sqrt{6}+\sqrt{2})/2\}$, then there are two different ways to realize $K_4$ as a $\{1,d\}$-graph, as shown in figure \ref{k4k5fig}. As we shall see shortly, this is very useful as in both cases we can represent $K_5\setminus e$, the complete graph on $5$ vertices with a missing edge, as a $\{1,d\}$-graph. Since in every $4$-coloring of $K_5\setminus e$ the endpoints of the missing edge must be colored the same, we can then employ the spindle technique to construct a $5$-chromatic $\{1,d\}$-graph. Details are presented below.

\begin{thm}\label{root3thm}
\begin{equation*}
\chi\left(\mathbb{E}^2, \left\{1, \sqrt{3}\right\}\right)=\chi\left(\mathbb{E}^2, \left\{1, 1/\sqrt{3}\right\}\right)\ge 5.
\end{equation*}
\end{thm}
\begin{proof}
Consider the following five points:
\begin{equation*}
\mathbf{1}\, (0,0),\,\,\mathbf{2}\,(2,0),\, \mathbf{3}\,(1/2, -\sqrt{3}/2),\, \mathbf{4}\, (1,0),\, \mathbf{5}\,(1/2,\sqrt{3}/2).
\end{equation*}
It is easy to see that with the exception of the pair $\{1,2\}$, every other two points are either $1$ or $\sqrt{3}$ apart. In other words, these points define a $\{1,\sqrt{3}\}$-graph
which is a complete $K_5$ with a missing edge, as shown in figure \ref{root3fig}.

\begin{figure}[htp]
\centering
\includegraphics[width=\linewidth]{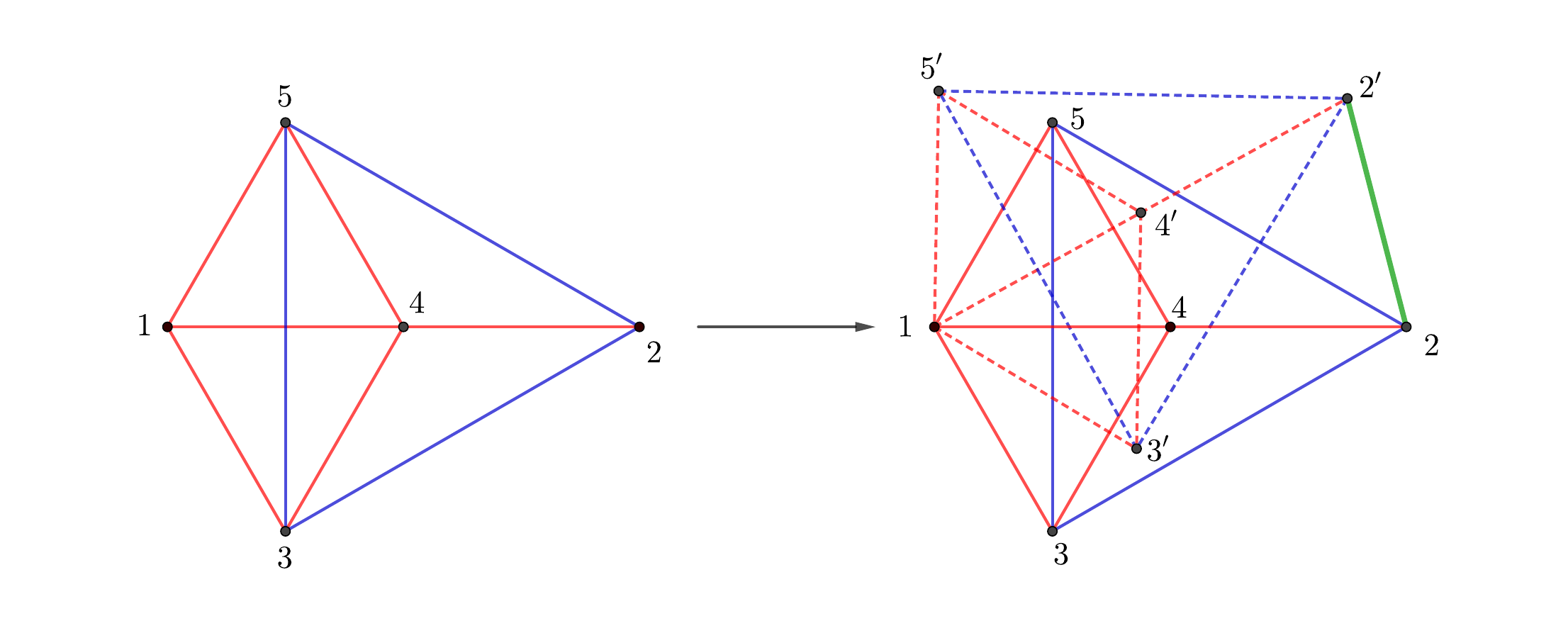}
\caption{ A $4$-colorable $K_5 \setminus e$ as a $\{1,\sqrt{3}\}$-graph converted to a $5$-chromatic  $\{1,\sqrt{3}\}$-graph via the spindle method}
\label{root3fig}
\end{figure}

It follows that for any $4$-coloring of this graph, vertices $1$ and $2$ must be assigned the same color. One can then apply the spindle technique by rotating the graph around vertex $1$ until the image of vertex $2$ is one unit apart from its original position. It is a simple exercise to find that the rotation angle can be chosen as $\arccos(7/8)$.
The coordinates of the four new vertices can be readily computed:
\begin{equation*}
\mathbf{2'}\,(7, \sqrt{15})/4,\,\mathbf{3'}\,(7+3\sqrt{5}, -7\sqrt{3}+\sqrt{15})/16,\,\mathbf{4'}\,(7, \sqrt{15})/8,\,\mathbf{5'}\,(7-3\sqrt{5}, 7\sqrt{3}+\sqrt{15})/16.
\end{equation*}
One obtains a $9$-vertex, $19$ edges, $\{1,\sqrt{3}\}$-graph, which requires $5$ colors.
\end{proof}

\newpage
One can use exactly the same approach to prove the following:

\begin{thm}
\begin{equation*}
\chi\left(\mathbb{E}^2, \left\{1, (\sqrt{6}+\sqrt{2})/2\right\}\right)=\chi\left(\mathbb{E}^2, \left\{1, (\sqrt{6}-\sqrt{2})/2\right\}\right)\ge 5.
\end{equation*}
\end{thm}

\begin{proof}
Consider the following five points:

\begin{equation*}
\mathbf{1}\, (0,0),\,\,\mathbf{2}\,(\sqrt{6},-\sqrt{2})/2,\, \mathbf{3}\,(-\sqrt{2}, -\sqrt{2})/2,\, \mathbf{4}\, (-\sqrt{2}+\sqrt{6},-\sqrt{2}-\sqrt{6})/4,\, \mathbf{5}\,(-\sqrt{2}+\sqrt{6},\sqrt{2}+\sqrt{6})/4.
\end{equation*}

It is easy to see that with the exception of the pair $\{1,2\}$, every other pair of points is either $1$ or $(\sqrt{6}+\sqrt{2})/2$ apart. In other words, these points define a $\{1,(\sqrt{6}+\sqrt{2})/2\}$-graph
which is a complete $K_5$ with a missing edge, as illustrated in figure \ref{root6fig}.

\begin{figure}[htp]
\centering
\includegraphics[width=\linewidth]{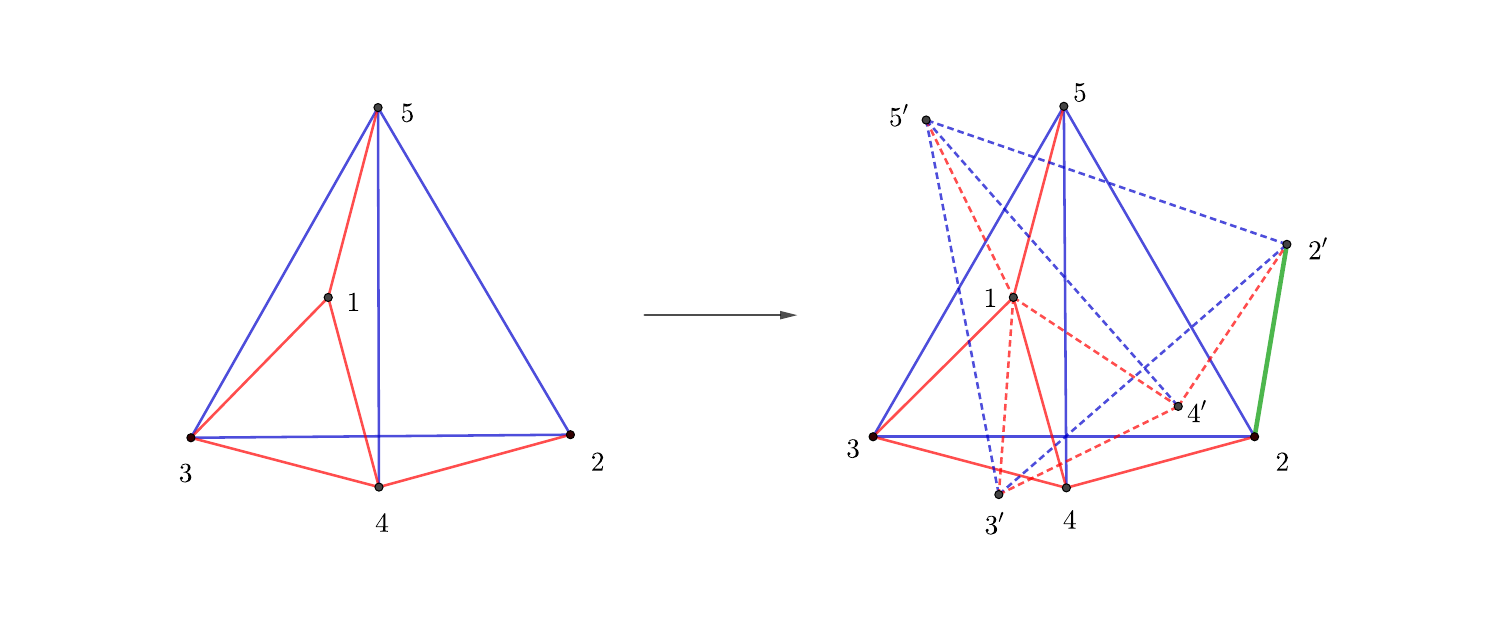}
\caption{ A $4$-colorable $K_5 \setminus e$ as a $\{1,(\sqrt{6}+\sqrt{2})/2\}$-graph converted to a $5$-chromatic  $\{1,(\sqrt{6}+\sqrt{2})/2\}$-graph via the spindle method}
\label{root6fig}
\end{figure}

It follows that for any $4$-coloring of this graph, vertices $1$ and $2$ must be assigned the same color. One can then apply the spindle technique by rotating the graph around vertex $1$ until the image of vertex $2$ is one unit apart from its original position. It is a simple exercise to find that the rotation angle can be chosen as $\arccos(3/4)$.
The coordinates of the new four vertices can be readily computed:

\begin{align*}
&\mathbf{2'}\,(3\sqrt{6}+\sqrt{14},-3\sqrt{2}+\sqrt{42})/8,\,\mathbf{3'}\, (-3\sqrt{2}+\sqrt{14},-3\sqrt{2}-\sqrt{14})/8,\\
&\mathbf{4'}\,(-3\sqrt{2}+3\sqrt{6}+\sqrt{14}+\sqrt{42}, -3\sqrt{2}-3\sqrt{6}-\sqrt{14}+\sqrt{42})/16,\\
&\mathbf{5'}\,(-3\sqrt{2}+3\sqrt{6}-\sqrt{14}-\sqrt{42},  \,3\sqrt{2}+3\sqrt{6}-\sqrt{14}+\sqrt{42})/16.
\end{align*}

One obtains a $9$-vertex, $19$ edges, $\{1,(\sqrt{6}+\sqrt{2})/2\}$-graph, which requires $5$ colors.
\end{proof}

\section{\bf The case $d=\sqrt{2}$.}

\begin{thm}\label{root2thm}
\begin{equation*}
\chi\left(\mathbb{E}^2, \left\{1, \sqrt{2}\right\}\right) \ge 5.
\end{equation*}
\end{thm}

As mentioned in the introduction, this result was proved by Katz, Krebs and Shaheen \cite{KKS}. They first proved the following:
\begin{thm}\label{kksthm}\cite{KKS}
Let $f:\mathbb{E}^2\rightarrow \mathbf{R}$ be a function such that $f(A)+f(B)+f(C)+f(D)=0$ holds whenever $ABCD$ is a unit square. Then $f(P)=0$ for all $P\in \mathbb{E}^2$.
\end{thm}
The proof is far from being easy. However, once established, it easily implies Theorem \ref{root2thm} as follows.

\begin{proof}\cite{KKS}
Assume that it is possible to color the plane with four colors (red, purple, green and blue), such that no two points distance $1$ or $\sqrt{2}$ apart are colored the same.
Observe that all four colors are needed for the vertices of a unit square. Define a function $f:\mathbb{E}^2\rightarrow \mathbf{R}$ by $f(P)=3$ if $P$ is colored purple, and $f(P)=-1$ otherwise. Then $f(A)+f(B)+f(C)+f(D)=0$ whenever $ABCD$ is a unit square. By Theorem \ref{kksthm}, $f$ must be the zero function, which contradicts the definition of $f$.
\end{proof}

The above proof is a beautiful example of an indirect argument since a $5$-chromatic $\{1,\sqrt{2}\}$-graph is not constructed explicitly. It is easy however to build such a graph, which we present below.
\begin{proof}
Consider the $\{1, \sqrt{2}\}$- graph whose vertices have the following coordinates:
\begin{align*}
&\mathbf{1}\,(1/2,1/2),\, \mathbf{2}\, (0,0),\, \mathbf{3}\, (1,0),\, \mathbf{4}\, (1,1),\,\mathbf{5}(0,1),\mathbf{6}\,(1/4-\sqrt{7}/4, 3/4-\sqrt{7}/4),\\
&\mathbf{7}\,(3/4+\sqrt{7}/4, 1/4+\sqrt{7}/4), \mathbf{8}\,(1/2, \sqrt{7}/2),\mathbf{9}\,(1/4-\sqrt{7}/4, 1/4+\sqrt{7}/4),\,\mathbf{10}\,(-\sqrt{7}/2, 1/2),\\
&\mathbf{11}(1+\sqrt{7}/2, 1/2),\,\mathbf{12}\,(3/4+\sqrt{7}/4, 3/4-\sqrt{7}/4),\,\mathbf{13}\,(1/2,1-\sqrt{7}/2).
\end{align*}
\begin{figure}[htp]
\centering
\includegraphics[width=0.8\linewidth]{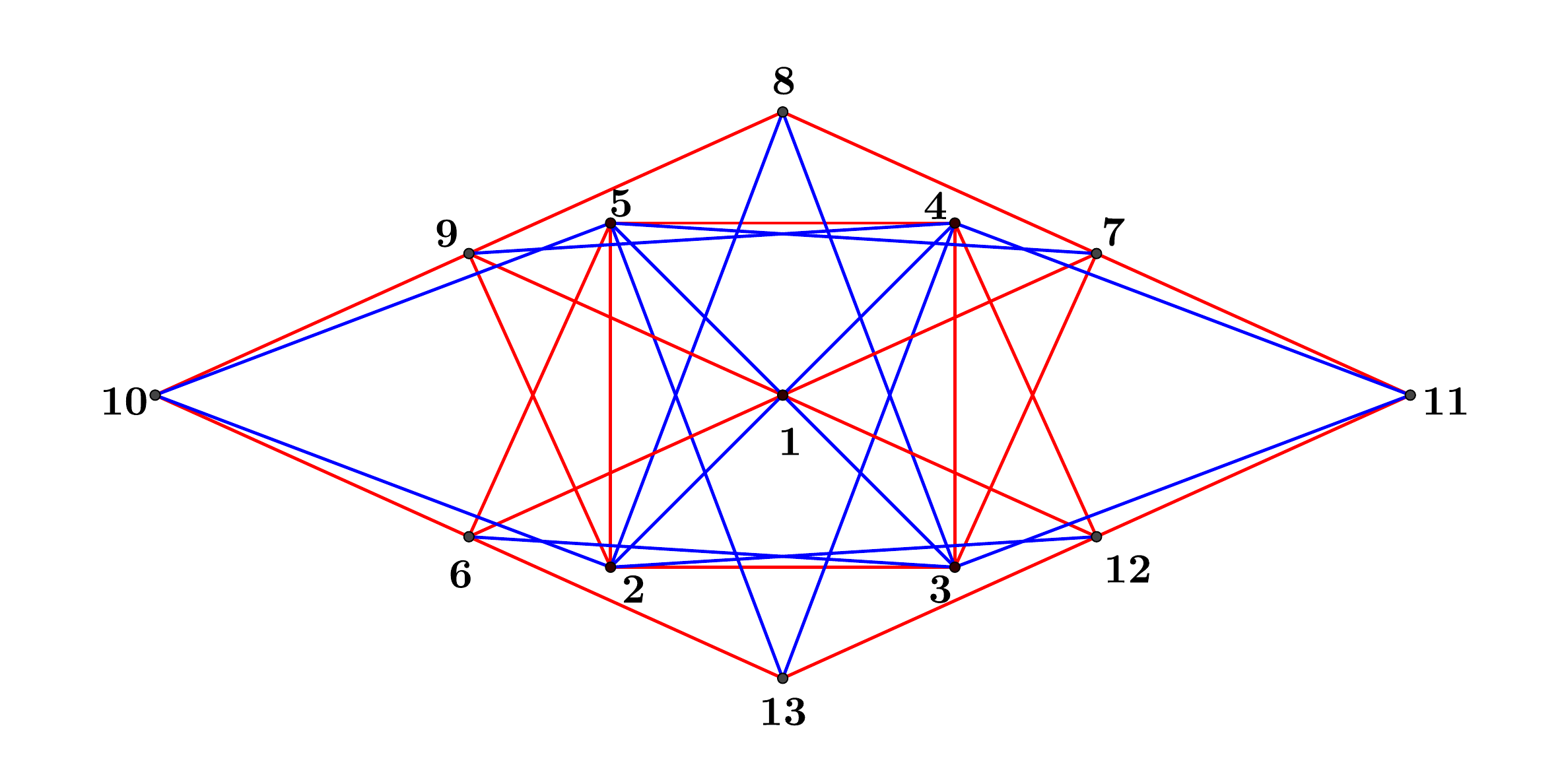}
\caption{Unit edges (in red) and $\sqrt{2}$-edges (in blue) of the $\{1, \sqrt{2}\}$-graph in Theorem \ref{root2thm}}
\label{root2fig}
\end{figure}

It is a simple matter to verify that the following $20$ pairs are unit-distance apart
\begin{align*}
&\{1, 6\}, \{1, 7\}, \{1, 9\}, \{1, 12\}, \{2, 3\}, \{2, 5\}, \{2, 9\}, \{3, 4\}, \{3, 7\}, \{4, 5\}, \{4, 12\}, \{5, 6\}, \{6, 10\}, \\
&\{6, 13\}, \{7, 8\}, \{7, 11\}, \{8, 9\}, \{9, 10\}, \{11, 12\}, \{12, 13\},
\end{align*}
while the following $14$ pairs are distance $\sqrt{2}$ apart:
\begin{equation*}
\{2, 4\}, \{2, 8\}, \{2, 10\}, \{2, 12\}, \{3, 5\}, \{3, 6\}, \{3, 8\}, \{3, 11\}, \{4, 9\}, \{4, 11\}, \{4, 13\}, \{5, 7\}, \{5, 10\}, \{5, 13\}.
\end{equation*}
We will show that this graph requires five colors. Let us assume that four colors suffice: \emph{purple, red, green and blue}. Since every two of the vertices $2, 3, 4$, and $5$ are
adjacent, one can assume that $c(2)=purple$, $c(3)=red$, $c(4)=green$, and $c(5)=blue$. By the symmetry of the graph one can assume that $c(1)=purple$, as well. But then $c(6)=c(7)=green$ since both vertices $6$ and $7$ are adjacent to $1$, $3$ and $5$. Since vertex $8$ is adjacent to $2, 3$, and $7$; then $c(8)=blue$. Similarly, vertex $9$ is adjacent to vertices $1, 4$, and $8$; hence $c(9)=red$. Vertex $10$ is adjacent to vertices $2, 5, 6,$ and $9$, which all have different colors. Hence a fifth color is needed. We reached the desired contradiction.
\end{proof}

\section{\bf An exotic looking distance} 

At this point we exhausted the possibilities that $K_4$, the complete graph on four vertices, can be represented as a $\{1,d\}$-graph. (see Theorem \ref{twodistance}).
Recall that we would like our graph to contain many small subgraphs which are hard to $4$-color. We already used $K_5 \setminus e$,
the complete graph on five vertices with a missing edge, in dealing with the cases $d=\sqrt{3}$ and $d=(\sqrt{6}+\sqrt{2})/2$.



There are no other useful subgraphs of order $5$ left to consider since $K_5 \setminus \{e, e'\}$, the $K_5$ with two missing edges can be $3$-colored if it is $K_4$-free.
Naturally, we look at subgraphs of order $6$. The question becomes: \emph{Are there any $4$-chromatic graphs of order $6$ which do not contain $K_4$ as a subgraph?}

The answer to the above question is affirmative. In fact, it is not difficult to show that there is exactly one graph with the desired properties: the wheel graph $W_6$.

There are several values of $d$ for which $W_6$ can be realized as a $\{1,d\}$-graph, including
\begin{equation}\label{other}
d=\sqrt{2+\sqrt{2}},\,\, d=\frac{1}{2}\sqrt{4+2\sqrt{2}},\,\, d=\frac{1}{2}\sqrt{8+2\sqrt{6}-2\sqrt{2}}.
\end{equation}

We will, however, consider a value which is somewhat more interesting. For the remainder of this section we will take
\begin{equation}\label{weird}
d= \frac{1}{2}\sqrt{3^{1/4}\cdot 2\sqrt{2}+2\sqrt{3}+2}.
\end{equation}

Note the $3^{1/4}$ which appears in the expression of $d$ above. The surprising fact is that despite the relatively complicated value of $d$, there are
$16$ different ways to embed $W_6$ as a $\{1,d\}$-graph - see figure \ref{16embeddings} in the appendix. This does not happen for any of the other distances listed in \eqref{other}.
It is therefore to be expected that it should be not too difficult to construct a $5$-chromatic $\{1,d\}$-graph. We show the details below.
\begin{lemma}\label{exoticlemma}
Consider the $\{1,d\}$- graph whose vertices have the following coordinates:
\begin{align*}
&\mathbf{1}\,(3^{1/4}\sqrt{2}-\sqrt{3}+1, 3^{3/4}\sqrt{2}+\sqrt{3}+1)/4, \,\, \mathbf{2}\, (-3^{1/4}\sqrt{2}+\sqrt{3}-1, 3^{3/4}\sqrt{2}+\sqrt{3}+1)/4,\\
&\mathbf{3}\,(3^{1/4}\sqrt{2}+\sqrt{3}+1, 3^{3/4}\sqrt{2}+\sqrt{3}-1)/4, \,\, \mathbf{4}\, (-3^{1/4}\sqrt{2}-\sqrt{3}-1, 3^{3/4}\sqrt{2}+\sqrt{3}-1)/4,\\
&\mathbf{5}\,(3^{1/4}\sqrt{2}+\sqrt{3}+1, 3^{3/4}\sqrt{2}+\sqrt{3}+3)/4, \,\, \mathbf{6}\, (-3^{1/4}\sqrt{2}-\sqrt{3}-1, 3^{3/4}\sqrt{2}+\sqrt{3}+3)/4,\\
&\mathbf{7}\,(3^{3/4}\sqrt{2}-3^{1/4}\sqrt{2}, 3^{3/4}\sqrt{2}-3^{1/4}\sqrt{2}+2\sqrt{3}-2)/4,\\
&\mathbf{8}\,(-3^{3/4}\sqrt{2}+3^{1/4}\sqrt{2}, 3^{3/4}\sqrt{2}-3^{1/4}\sqrt{2}+2\sqrt{3}-2)/4,\\
&\mathbf{9}\, (1,0), \, \mathbf{10}\,(-1,0),\, \mathbf{11}\, (1, \sqrt{3})/2,\, \mathbf{12}\, (-1,\sqrt{3})/2,\, \mathbf{13}\,(0,0).
\end{align*}
Then in every $4$-coloring of this graph, vertices $1$ and $2$ must be assigned the same color.
\end{lemma}
\vspace{-0.5cm}
\begin{figure}[htp]
\centering
\includegraphics[width=0.7\linewidth]{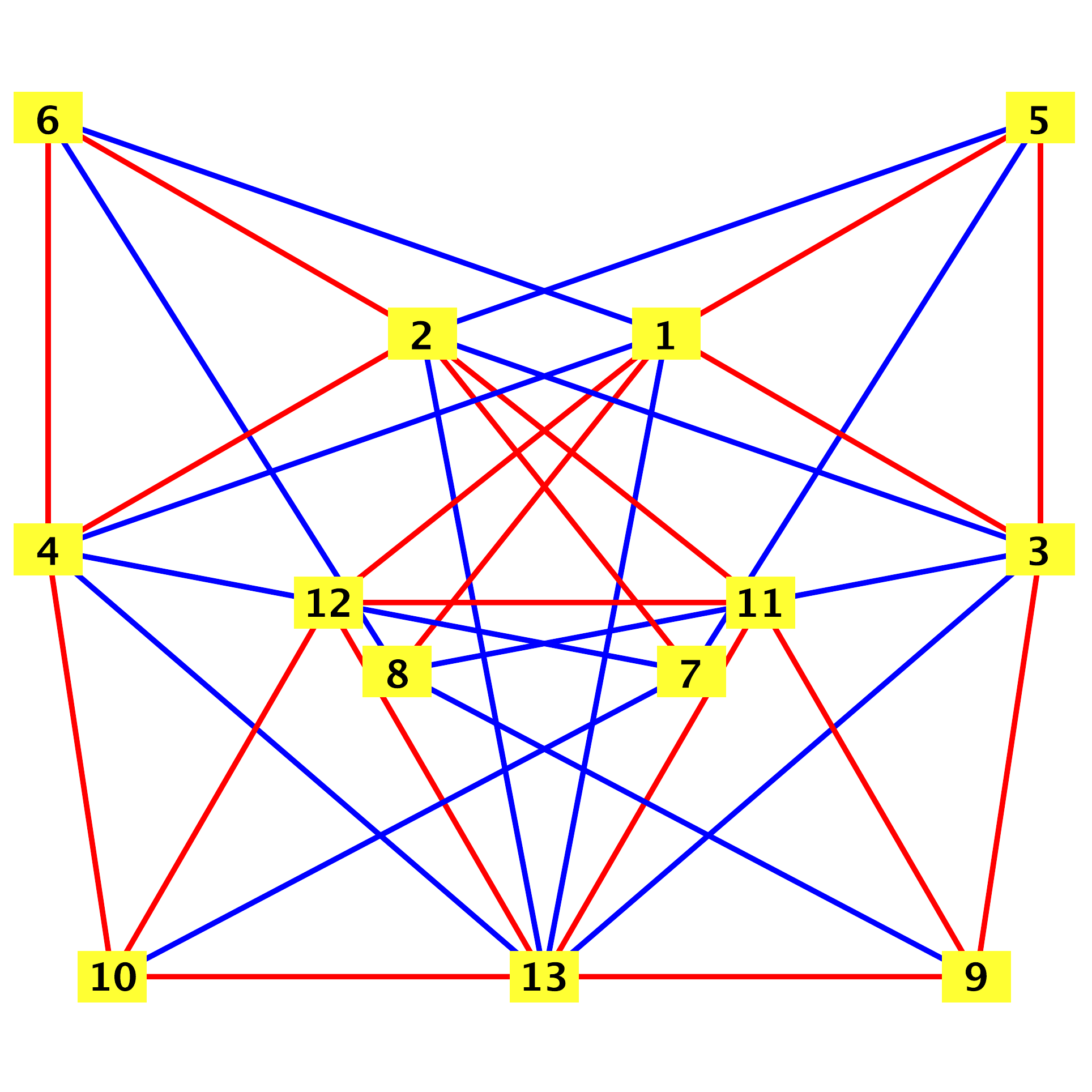}
\vspace{-0.3cm}
\caption{Unit edges (in red) and $d$-edges (in blue) for the graph in Lemma \ref{exoticlemma}}
\label{exotic}
\end{figure}
\begin{proof}
It is easy to verify that the thirteen points listed above determine $19$ unit edges:
\begin{align*}
&\{1, 3\}, \{1, 5\}, \{1, 8\}, \{1, 12\}, \{2, 4\}, \{2, 6\}, \{2, 7\}, \{2, 11\}, \{3, 5\}, \{3, 9\}, \{4, 6\}, \{4, 10\}, \{9, 11\}, \\
&\{9, 13\}, \{10, 12\}, \{10, 13\}, \{11, 12\}, \{11, 13\}, \{12, 13\},
\end{align*}
and $14$ edges of length $d$:
$$\{1, 4\}, \{1, 6\}, \{1, 13\}, \{2, 3\}, \{2, 5\}, \{2, 13\}, \{3, 8\}, \{3, 13\}, \{4, 7\}, \{4, 13\}, \{5, 7\}, \{6, 8\}, \{7, 10\}, \{8, 9\}.$$
Assume by contradiction that there exists a $4$-coloring $c:\{1,2,\ldots,12,13\}\rightarrow \{red, purple, green, blue\}$ such that $c(1)\neq c(2)$.
Without loss of generality, assume that $c(1)=purple$ and $c(2)=red$. Since vertices $3$ and $13$ are both adjacent to vertices
$1$ and $2$ and adjacent to each other, we can assume that $c(3)=green$ and $c(13)=blue$. With this setup, we can now show that the colors of the remaining vertices are uniquely determined.

Vertex $4$ is adjacent to $1$, $2$, and $13$, hence $c(4)=green$.
Vertex $5$ is adjacent to $1$, $2$, and $3$, hence $c(5)=blue$.
Vertex $6$ is adjacent to $1$, $2$, and $4$, hence $c(6)=blue$.
Vertex $7$ is adjacent to $2$, $4$, and $5$, hence $c(7)=purple$.
Vertex $8$ is adjacent to $1$, $3$, and $6$, hence $c(8)=red$.
Vertex $9$ is adjacent to $3$, $8$, and $13$, hence $c(9)=purple$.
Vertex $10$ is adjacent to $4$, $7$, and $13$, hence $c(10)=red$.
Vertex $11$ is adjacent to $2$, $9$, and $13$, hence $c(11)=green$.
Vertex $12$ is adjacent to $1$, $10$, and $13$, hence $c(12)=green$.

However, $c(11)=c(12)=green$ is impossible since vertices $11$ and $12$ are adjacent.
We obtained the desired contradiction.
\end{proof}
We are now in position to prove the following
\begin{thm}
\begin{equation*}
\chi\left(\mathbb{E}^2, \left\{1, \frac{1}{2}\sqrt{3^{1/4}\cdot 2\sqrt{2}+2\sqrt{3}+2}\right\}\right) \ge 5.
\end{equation*}
\end{thm}
\begin{proof}
Consider the graph defined in Lemma \ref{exoticlemma}. Any $4$-coloring must assign vertices $1$ and $2$ the same color. Moreover, the distance between these two vertices
is $(3^{1/4}\sqrt{2}-\sqrt{3}+1)/2 = 0.564\ldots >1/2$. Hence, we can apply the spindle technique to this graph by rotating about vertex $1$ until the image of vertex $2$ is at distance $1$ from its original position. The resulting $\{1,d\}$-graph has $25$ vertices, $67$ edges, and it is $5$-chromatic.
\end{proof}

\section{\bf A special subset of $\mathbb{E}^2$}

For the rest of the paper we focus on graphs whose vertices belong to the following subset of $\mathbb{E}^2$:
\begin{equation}\label{lambda}
\Lambda:=\left\{\left(\frac{a\sqrt{3}}{12}+\frac{b\sqrt{11}}{12}, \frac{c}{12}+\frac{d\sqrt{33}}{12}\right)\,:\, a, b, c, d, \text{integers}\right\}.
\end{equation}
To simplify the presentation, in the sequel we will use the following notation:
\begin{equation}\label{convention}
[a,b,c,d]:=\left(\frac{a\sqrt{3}}{12}+\frac{b\sqrt{11}}{12}, \frac{c}{12}+\frac{d\sqrt{33}}{12}\right).
\end{equation}
One reason for this choice is that Mosers' $7$-vertex $4$-chromatic unit
distance graph, shown in figure \ref{spindlefig}, has such an embedding.  With
the convention in \eqref{convention}, the coordinates of the vertices can be
written
\begin{equation}
[0, 0, 0, 0], [0, 0, 12, 0], [6, 0, 6, 0], [6, 0, 18, 0], [0, 2, 10, 0], [5, 1, 5, -1], [5, 3, 15, -1].
\end{equation}

\newpage
\section{\bf Two more distances, and a slightly different approach}
The next two results are easy to prove; however, we believe the method of proof is significant.
\begin{lemma}\label{smart1}
If $\chi\left(\mathbb{E}^2,\left\{1, \sqrt{3/2+\sqrt{33}/6}\right\}\right)=4$ then $\chi\left(\mathbb{E}^2, \left\{1, \sqrt{3/2+\sqrt{33}/6},1/\sqrt{3}\right\}\right)=4$.
\end{lemma}
\begin{proof}
Suppose that $\chi\left(\mathbb{E}^2,\left\{1,
\sqrt{3/2+\sqrt{33}/6}\right\}\right)=4$; that is, there exists a $4$-coloring
of the plane, $c:\mathbb{E}^2\rightarrow \{red, purple, green, blue\}$ such
that no two points which are either distance $1$ or distance
$\sqrt{3/2+\sqrt{33}/6}$ apart are assigned the same color. We will show that
for this particular coloring, no two points distance $1/\sqrt{3}$ apart can be
colored identically either.
Indeed, for the sake of reaching a contradiction, assume that two such points
exist, and denote them by $A$ and $B$.  Choose a system of coordinates such
that $A=(0,0)=[0,0,0,0]$ and $B=(1/\sqrt{3},0)=[4,0,0,0]$. Here we used
notation \eqref{convention} for the coordinates.  Consider next the
$\{1,\sqrt{3/2+\sqrt{33}/6}\}$-graph induced by $A$, $B$, and the  following
seven additional vertices:
\begin{equation*}
\mathbf{1}\,[1, 3, 3, -1], \mathbf{2}\,[0, 0, 12, 0], \mathbf{3}\,[-1, -3, 9, 1], \mathbf{4}\,[-1, 3, 3, 1], \mathbf{5}\,[-2, 0, 6, 0], \mathbf{6}\,[0, 0, 6, 2], \mathbf{7}\,[5, 3, 9, 1].
\end{equation*}

Assume that $c(A)=c(B)=blue$, as shown in figure \ref{3over2fig}. Note that vertex $A$ is adjacent to vertices $1$, $2$, $3$, and $4$, while vertex $B$ is adjacent to vertices $5$, $6$ and $7$. Hence, none of vertices $1$ through $7$ can be colored blue; that is, we can only use the remaining three colors for these vertices. However, it is easy to see that the subgraph induced by these seven vertices cannot be $3$-colored.

\vspace{-0.1cm}
\begin{figure}[htp]
\centering
\includegraphics[width=0.58\linewidth]{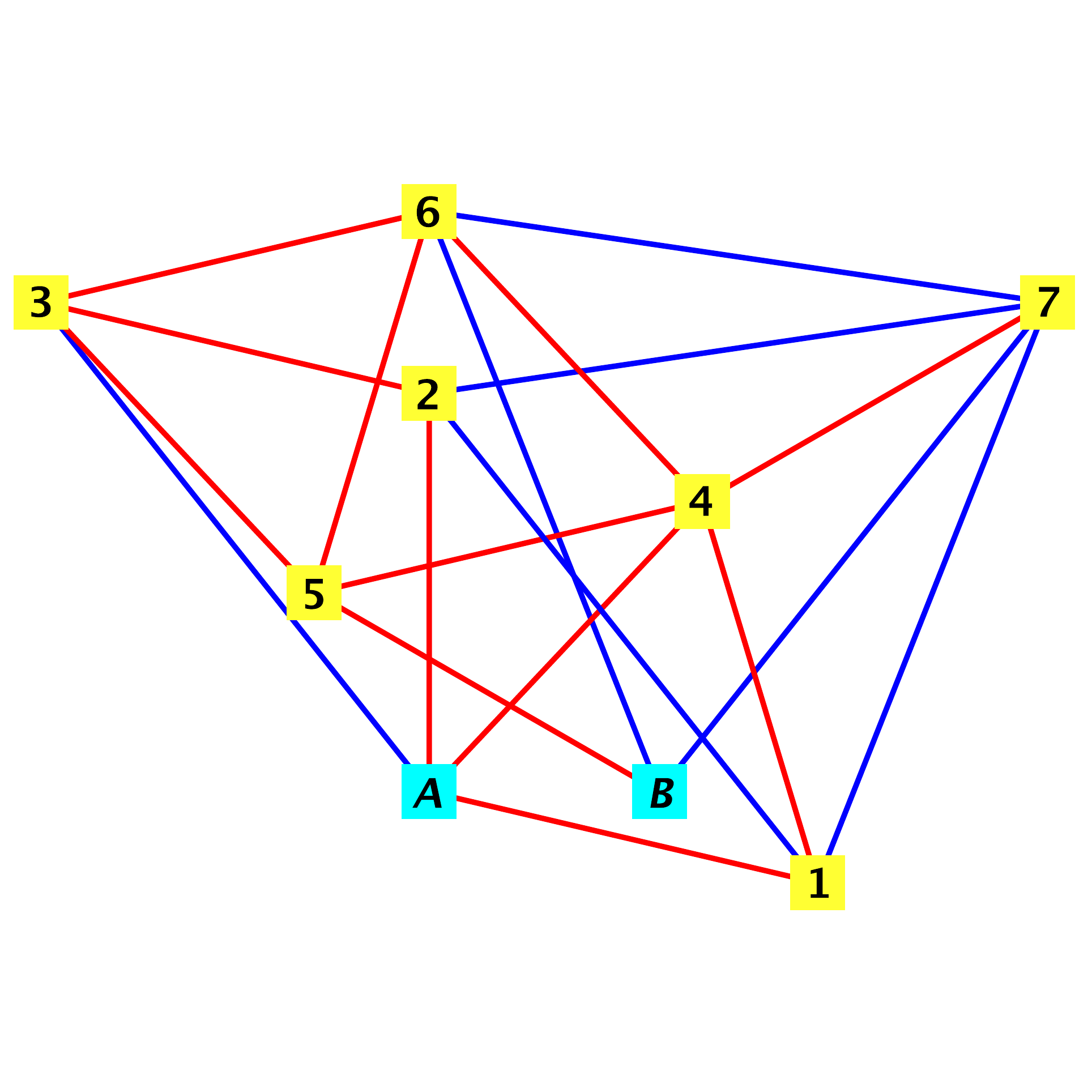}
\vspace{-1cm}
\caption{A $\left\{1,\sqrt{3/2+\sqrt{33}/6}\right\}$-graph. Unit edges are in red, $\sqrt{3/2+\sqrt{33}/6}$ edges are in blue. The distance between points $A$ and $B$ is $1/\sqrt{3}$.}
\label{3over2fig}
\end{figure}

Indeed, since vertices $3$, $5$, and $6$ form a triangle, one can assume that
$c(3)=red$, $c(5)=green$, and $c(6)=purple$. But this choice implies in order
that $c(4)=red$, $c(7)=green$, and $c(2)=purple$. Finally, vertex $1$ is
adjacent to $2$, $4$, and $7$, and these vertices are already colored with
different colors. Hence, the subgraph induced by vertices $1$ through $7$
cannot be $3$-colored.  We have obtained the desired contradiction.
\end{proof}

The above result combined with Theorem \ref{root3thm} immediately implies the following:
\begin{thm}
$\chi\left(\mathbb{E}^2,\left\{1, \sqrt{3/2+\sqrt{33}/6}\right\}\right)\ge 5$.
\end{thm}
\begin{proof}
Suppose the opposite is true.

Then by Lemma \ref{smart1} it follows that $\chi\left(\mathbb{E}^2, \left\{1, \sqrt{3/2+\sqrt{33}/6},1/\sqrt{3}\right\}\right)=4$, with eventually gives that $\chi\left(\mathbb{E}^2, \left\{1, 1/\sqrt{3}\right\}\right)=4$. But this contradicts the result in Theorem \ref{root3thm}.
\end{proof}

{\bf Observation.}
One can construct an explicit $5$-chromatic $\{1,\sqrt{3/2+\sqrt{33}/6}\}$-graph as follows:

Start with the $9$-vertex graph from Theorem \ref{root3thm} scaled down by a factor of $\sqrt{3}$. This graph has $13$ edges of length $1/\sqrt{3}$ and $6$ edges of length $1$, and it is a $5$-chromatic $\{1, 1/\sqrt{3}\}$-graph.

Next, for each edge $AB$ of length $1/\sqrt{3}$, add to this graph the seven vertices of the graph constructed in Lemma \ref{smart1}.
The resulting graph has $9+13\cdot 7=100$ vertices, and it is a $5$-chromatic $\{1,\sqrt{3/2+\sqrt{33}/6}\}$-graph.
Most likely, much smaller graphs with this property do exist.  The advantage of our proof is that we do not have to deal with the $100$-vertex graph directly; instead
a lower bound for the chromatic number of this graph follows by considering two much smaller graphs, both of order $9$.

Using the same technique as in Lemma \ref{smart1} we can prove the following

\begin{lemma}\label{smart2}
If $\chi\left(\mathbb{E}^2,\left\{1, \sqrt{5/3}\right\}\right)=4$ then $\chi\left(\mathbb{E}^2, \left\{1, \sqrt{5/3},1/\sqrt{3}\right\}\right)=4$.
\end{lemma}
\begin{proof}
Suppose that $\chi\left(\mathbb{E}^2,\left\{1, \sqrt{5/3}\right\}\right)=4$; that is, there exists a $4$-coloring of the plane $c:\mathbb{E}^2\rightarrow \{red, purple, green, blue\}$ such that no two points which are either distance $1$ or distance $\sqrt{5/3}$ apart are assigned the same color. We will show that for this particular coloring, no two points distance $1/\sqrt{3}$ apart can be colored identically either.

In order to a contradiction, assume that two such points exist and denote them by $A$ and $B$.  Choose a system of coordinates such that $A=(0,0)=[0,0,0,0]$ and $B=(1/\sqrt{3},0)=[4,0,0,0]$. Again, we are using notation \eqref{convention} for the coordinates.

\begin{figure}[htp]
\centering
\includegraphics[width=0.8\linewidth]{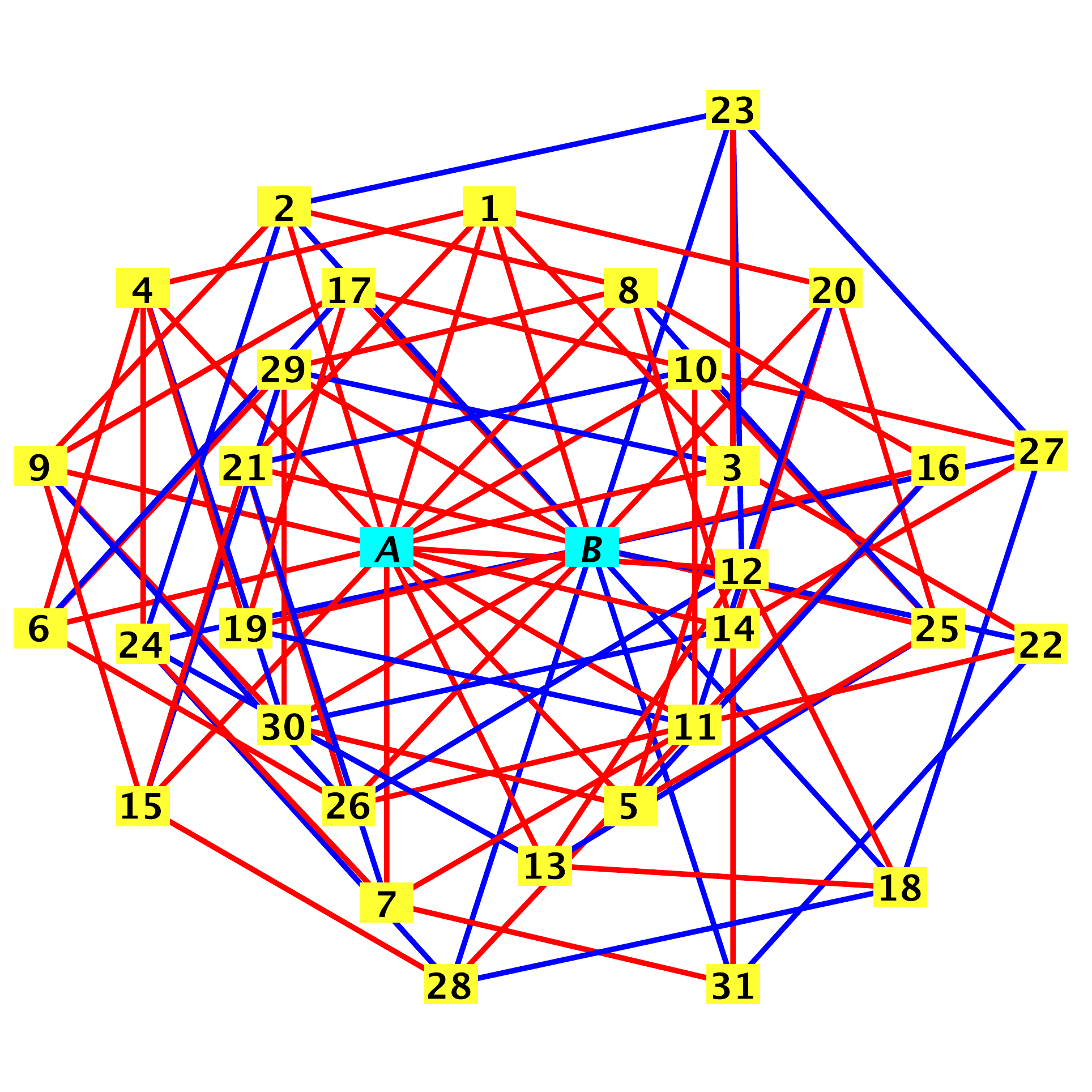}
\caption{A $33$-vertex $\left\{1,\sqrt{5/3}\right\}$-graph. Unit edges are in red; $\sqrt{5/3}$ edges are in blue. The distance between points $A$ and $B$ is $1/\sqrt{3}$.}
\label{5over3fig}
\end{figure}

Consider next the $\{1,\sqrt{5/3}\}$-graph induced by $A$, $B$, and the  following $31$ additional vertices:

\begin{align*}
&[2, 0, 0, 2], [-2, 0, 0, 2], [1, 3, -3, 1], [1, -3, 3, 1], [-1, 3, -3, -1], [-1, -3, 3, -1], [0, 0, -12, 0], [-1, 3, 3, 1],\\
&[-1, -3, -3, 1], [6, 0, 6, 0], [6, 0, -6, 0], [5, 1, 5, -1], [5, -1, -5, -1], [1, 3, 3, -1], [1, -3, -3, -1], [5, 3, -3, 1],\\
&[5, -3, 3, 1], [10, 0, 0, -2], [3, -3, 3, -1], [3, 3, 3, 1], [3, -3, -3, 1], [7, 3, -9, 1], [1, 3, 9, 1], [1, -3, -9, 1],\\
&[5, 3, 3, -1], [5, -3, -3, -1], [7, 3, 9, -1], [7, -3, -9, -1], [-2, 0, 6, 0], [-2, 0, -6, 0], [1, 3, -9, -1].
\end{align*}

Assume that $c(A)=c(B)=blue$, as shown in figure \ref{5over3fig}. It can be
verified that vertex $A$ is adjacent to vertices $1$ through $15$, while vertex
$B$ is adjacent to vertices $16$ through $31$ (and to vertices $1$ and $2$ as
well). Hence, none of vertices $1$ through $31$ can be colored blue.  So
we can only use the remaining three colors for these vertices. However, it can
be checked that the subgraph induced by these $31$ vertices, regarded as a
$\{1,\sqrt{5/3}\}$-graph, cannot be $3$-colored. We used both Maple and Sage
to verify this assertion.
\end{proof}

Again, we combine the above Lemma with Theorem \ref{root3thm} to obtain the following
\begin{thm}
$\chi\left(\mathbb{E}^2,\left\{1, \sqrt{5/3}\right\}\right)\ge 5$.
\end{thm}
\begin{proof}
Suppose the opposite is true.

Then by Lemma \ref{smart2} it follows that $\chi\left(\mathbb{E}^2, \left\{1, \sqrt{5/3},1/\sqrt{3}\right\}\right)=4$, which eventually gives that $\chi\left(\mathbb{E}^2, \left\{1, 1/\sqrt{3}\right\}\right)=4$. This contradicts the result in Theorem \ref{root3thm}.
\end{proof}

\newpage
\section{\bf Another two distances}

One may argue that quite a few of the values of $d$ we studied so far are rather complicated. In this section we are trying to rectify this situation by considering
two simple values: $d=2$ and $d=2/\sqrt{3}$.
\begin{thm}
\begin{equation}
\chi(\mathbb{E}^2,\{1,2\})\ge 5.
\end{equation}
\end{thm}
\begin{proof}
Consider the $26$-vertex $\{1,2\}$-graph whose vertices are given by the following coordinates:
\begin{align*}
&[-2, 0, 0, 2], [2, 0, 0, 2], [0, 0, 0, 0], [0, 0, 0, 4], [0, 0, -6, 2], [0, 0, 6, 2], [-1, -3, 3, 3], [1, 3, 3, 3], [-3, -3, 3, 1],\\
&[3, 3, 3, 1], [-1, -3, -3, 1], [1, 3, -3, 1], [-4, 0, 0, 0], [4, 0, 0, 0], [3, -3, -3, 1], [-3, 3, -3, 1], [1, -3, -3, 3],\\
&[-1, 3, -3, 3], [1, -3, 3, 1], [-1, 3, 3, 1], [-2, 0, 6, 0], [2, 0, 6, 0], [-2, 0, -6, 0], [2, 0, -6, 0], [0, -6, 0, 2], [0, 6, 0, 2].
\end{align*}
\vspace{-0.5cm}
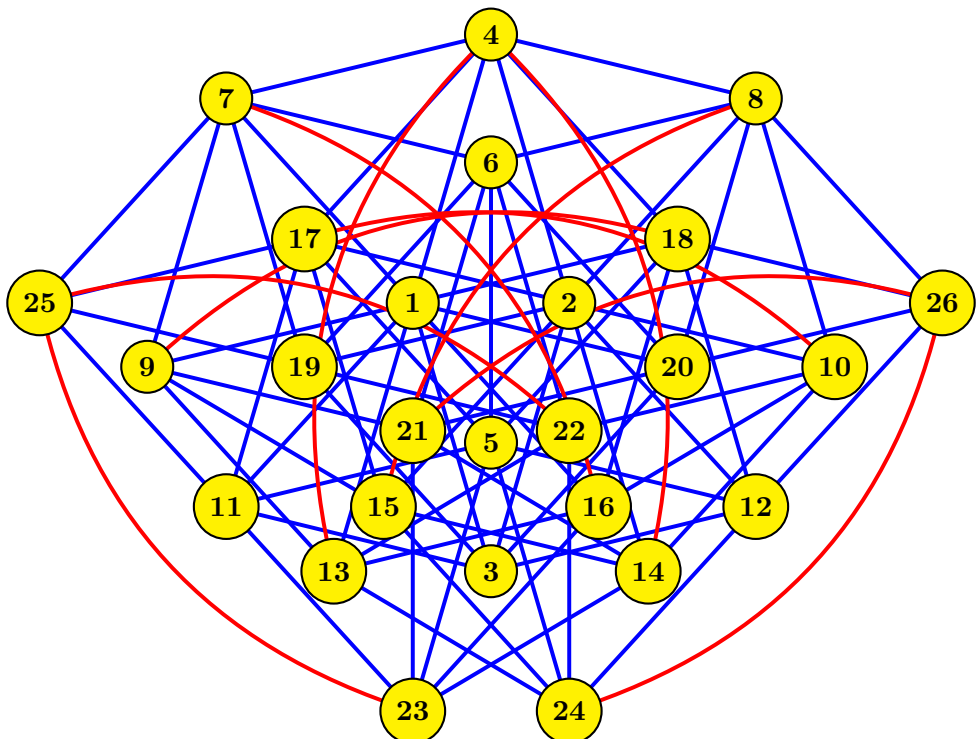
\begin{figure}[htp]
\centering
\begin{tikzpicture}[line width=1.5pt,state/.style={circle, draw, minimum size=0.65cm}]
\tikzstyle{fontbf} = [draw,text centered,font=\bf]
\tikzstyle{every node}=[draw=black,thick, shape=circle]

\coordinate (AA) at ( 4.96, 5.43);
\coordinate (AB) at ( 7.04, 5.43);
\coordinate (AC) at ( 6.00, 1.86);
\coordinate (AD) at ( 6.00, 9.00);
\coordinate (AE) at ( 6.00, 3.57);
\coordinate (AF) at ( 6.00, 7.30);
\coordinate (AG) at ( 2.48, 8.15);
\coordinate (AH) at ( 9.52, 8.15);
\coordinate (AI) at ( 1.43, 4.58);
\coordinate (AJ) at (10.57, 4.58);
\coordinate (AK) at ( 2.48, 2.72);
\coordinate (AL) at ( 9.52, 2.72);
\coordinate (AM) at ( 3.91, 1.86);
\coordinate (AN) at ( 8.09, 1.86);
\coordinate (AO) at ( 4.57, 2.72);
\coordinate (AP) at ( 7.43, 2.72);
\coordinate (AQ) at ( 3.52, 6.28);
\coordinate (AR) at ( 8.48, 6.28);
\coordinate (AS) at ( 3.52, 4.58);
\coordinate (AT) at ( 8.48, 4.58);
\coordinate (AU) at ( 4.96, 3.73);
\coordinate (AV) at ( 7.04, 3.73);
\coordinate (AW) at ( 4.96, 0.00);
\coordinate (AX) at ( 7.04, 0.00);
\coordinate (AY) at ( 0.00, 5.43);
\coordinate (AZ) at (12.00, 5.43);
\draw[blue] (AA) to (AC);
\draw[blue] (AA) to (AD);
\draw[blue] (AA) to (AG);
\draw[blue] (AA) to (AI);
\draw[blue] (AA) to (AK);
\draw[blue] (AA) to (AM);
\draw[blue] (AA) to (AP);
\draw[blue] (AA) to (AR);
\draw[blue] (AA) to (AT);
\draw[blue] (AB) to (AC);
\draw[blue] (AB) to (AD);
\draw[blue] (AB) to (AH);
\draw[blue] (AB) to (AJ);
\draw[blue] (AB) to (AL);
\draw[blue] (AB) to (AN);
\draw[blue] (AB) to (AO);
\draw[blue] (AB) to (AQ);
\draw[blue] (AB) to (AS);
\draw[blue] (AC) to (AK);
\draw[blue] (AC) to (AL);
\draw[blue] (AC) to (AS);
\draw[blue] (AC) to (AT);
\draw[blue] (AD) to (AG);
\draw[blue] (AD) to (AH);
\draw[blue] (AD) to (AQ);
\draw[blue] (AD) to (AR);
\draw[blue] (AE) to (AF);
\draw[blue] (AE) to (AK);
\draw[blue] (AE) to (AL);
\draw[blue] (AE) to (AQ);
\draw[blue] (AE) to (AR);
\draw[blue] (AE) to (AW);
\draw[blue] (AE) to (AX);
\draw[blue] (AF) to (AG);
\draw[blue] (AF) to (AH);
\draw[blue] (AF) to (AS);
\draw[blue] (AF) to (AT);
\draw[blue] (AF) to (AU);
\draw[blue] (AF) to (AV);
\draw[blue] (AG) to (AI);
\draw[blue] (AG) to (AS);
\draw[blue] (AG) to (AY);
\draw[blue] (AH) to (AJ);
\draw[blue] (AH) to (AT);
\draw[blue] (AH) to (AZ);
\draw[blue] (AI) to (AM);
\draw[blue] (AI) to (AO);
\draw[blue] (AI) to (AU);
\draw[blue] (AJ) to (AN);
\draw[blue] (AJ) to (AP);
\draw[blue] (AJ) to (AV);
\draw[blue] (AK) to (AQ);
\draw[blue] (AK) to (AW);
\draw[blue] (AK) to (AY);
\draw[blue] (AL) to (AR);
\draw[blue] (AL) to (AX);
\draw[blue] (AL) to (AZ);
\draw[blue] (AM) to (AP);
\draw[blue] (AM) to (AV);
\draw[blue] (AM) to (AX);
\draw[blue] (AN) to (AO);
\draw[blue] (AN) to (AU);
\draw[blue] (AN) to (AW);
\draw[blue] (AO) to (AQ);
\draw[blue] (AO) to (AX);
\draw[blue] (AP) to (AR);
\draw[blue] (AP) to (AW);
\draw[blue] (AQ) to (AY);
\draw[blue] (AR) to (AZ);
\draw[blue] (AS) to (AV);
\draw[blue] (AS) to (AY);
\draw[blue] (AT) to (AU);
\draw[blue] (AT) to (AZ);
\draw[blue] (AU) to (AW);
\draw[blue] (AV) to (AX);
\draw[red] (AD) to[bend right] (AM);
\draw[red] (AD) to[bend left] (AN);
\draw[red] (AG) to[bend left] (AP);
\draw[red] (AH) to[bend right] (AO);
\draw[red] (AI) to[bend left] (AR);
\draw[red] (AJ) to[bend right] (AQ);
\draw[red] (AU) to[bend left] (AZ);
\draw[red] (AV) to[bend right] (AY);
\draw[red] (AW) to[bend left] (AY);
\draw[red] (AX) to[bend right] (AZ);
\node[fill=yellow] at (AA) {\textbf{1}};
\node[fill=yellow] at (AB) {\textbf{2}};
\node[fill=yellow] at (AC) {\textbf{3}};
\node[fill=yellow] at (AD) {\textbf{4}};
\node[fill=yellow] at (AE) {\textbf{5}};
\node[fill=yellow] at (AF) {\textbf{6}};
\node[fill=yellow] at (AG) {\textbf{7}};
\node[fill=yellow] at (AH) {\textbf{8}};
\node[fill=yellow] at (AI) {\textbf{9}};
\node[fill=yellow] at (AJ) {\textbf{10}};
\node[fill=yellow] at (AK) {\textbf{11}};
\node[fill=yellow] at (AL) {\textbf{12}};
\node[fill=yellow] at (AM) {\textbf{13}};
\node[fill=yellow] at (AN) {\textbf{14}};
\node[fill=yellow] at (AO) {\textbf{15}};
\node[fill=yellow] at (AP) {\textbf{16}};
\node[fill=yellow] at (AQ) {\textbf{17}};
\node[fill=yellow] at (AR) {\textbf{18}};
\node[fill=yellow] at (AS) {\textbf{19}};
\node[fill=yellow] at (AT) {\textbf{20}};
\node[fill=yellow] at (AU) {\textbf{21}};
\node[fill=yellow] at (AV) {\textbf{22}};
\node[fill=yellow] at (AW) {\textbf{23}};
\node[fill=yellow] at (AX) {\textbf{24}};
\node[fill=yellow] at (AY) {\textbf{25}};
\node[fill=yellow] at (AZ) {\textbf{26}};

\end{tikzpicture}

\caption{A $26$-vertex $\left\{1,2\right\}$-graph with chromatic number $5$}
\label{d2}
\end{figure}

It can be verified that this graph has $75$ unit edges, and $10$ edges of length $2$.
The graph is shown in figure \ref{d2},
and it can be checked that it is $5$-chromatic. Again, we used both Maple and Sage to verify this.
Given the relatively small order of this graph, a computer-free proof is certainly possible.
Since the midpoint of any edge of length $2$ is also a vertex of the graph, these long edges are shown with curved line
segments.

The surprising fact is that a relatively small number of long edges (only $10$ of them), is sufficient to raise the chromatic number from $4$ to $5$.
\end{proof}

\begin{thm}\label{2overroot3thm}
\begin{equation}
\chi(\mathbb{E}^2,\{1,2/\sqrt{3}\})\ge 5.
\end{equation}
\end{thm}
\begin{proof}
Consider the $103$-vertex $\{1,2/\sqrt{3}\}$-graph whose vertices are given by the following coordinates:
{\tiny
\begin{align*}
&[0, 0, 0, 0], [6, 0, 6, 0], [0, 0, 12, 0], [-6, 0, 6, 0], [-6, 0, -6, 0], [0, 0, -12, 0], [6, 0, -6, 0], [0, -2, -10, 0], [0, 2, -10, 0], [0, -2, 10, 0], [0, 2, 10, 0], \\
&[-2, 0, 0, -2], [2, 0, 0, -2],[-2, 0, 0, 2], [2, 0, 0, 2], [-5, 1, -5, -1], [5, -1, -5, -1], [-5, -1, -5, 1], [5, 1, -5, 1], [-1, 3, -3, -1], [1, -3, -3, -1],\\
& [-1, -3, -3, 1], [1, 3, -3, 1], [-1, -3, 3, -1], [1, 3, 3, -1], [-1, 3, 3, 1], [1, -3, 3, 1], [-5, -1, 5, -1], [5, 1, 5, -1], [-5, 1, 5, 1], [5, -1, 5, 1], \\
&[0, -6, -6, 0], [0, 6, -6, 0], [0, -6, 6, 0], [0, 6, 6, 0], [-3, 3, -3, -3], [3, -3, -3, -3], [-3, -3, -3, 3], [3, 3, -3, 3], [-3, -3, 3, -3], [3, 3, 3, -3], \\
&[-3, 3, 3, 3], [3, -3, 3, 3], [-4, 0, 0, 0], [4, 0, 0, 0], [-2, 0, -6, 0], [2, 0, -6, 0], [-2, 0, 6, 0], [2, 0, 6, 0], [0, -2, -2, 0], [0, 2, -2, 0], [0, -2, 2, 0],\\
&[0, 2, 2, 0], [1, -1, -1, -1], [-1, 1, -1, -1], [-1, -1, -1, 1], [1, 1, -1, 1], [-1, -1, 1, -1], [1, 1, 1, -1], [1, -1, 1, 1], [-1, 1, 1, 1], [8, 0, 0, 0],\\
& [4, 0, 12, 0],[-4, 0, 12, 0], [-8, 0, 0, 0], [-4, 0, -12, 0], [4, 0, -12, 0], [0, -4, -4, 0], [0, 4, -4, 0], [0, -4, 4, 0], [0, 4, 4, 0], [-2, 2, -2, -2],\\
&[2, -2, -2, -2],[-2, -2, -2, 2], [2, 2, -2, 2], [-2, -2, 2, -2], [2, 2, 2, -2],[-2, 2, 2, 2], [2, -2, 2, 2], [-4, 2, -2, 0], [4, -2, -2, 0], [-4, 2, 2, 0], \\
&[4, -2, 2, 0], [1, -1, -7, 1],[-1, 1, -7, 1], [-3, 1, -5, 1], [3, -1, -5, 1], [-3, 1, 5, -1],[3, -1, 5, -1], [1, -1, 7, -1], [-1, 1, 7, -1], [-4, -2, -2, 0],\\
&[4, 2, -2, 0], [-4, -2, 2, 0], [4, 2, 2, 0],[-1, -1, -7, -1], [1, 1, -7, -1], [-3, -1, -5, -1], [3, 1, -5, -1], [-3, -1, 5, 1], [3, 1, 5, 1], [-1, -1, 7, 1], [1, 1, 7, 1].
\end{align*}}
\vspace{-1cm}

\begin{figure}
\begin{tikzpicture}[line width=1pt]
\tikzstyle{every node}=[draw=black,fill=yellow,
  shape=circle,minimum height=0.25cm,inner sep=1];

\input{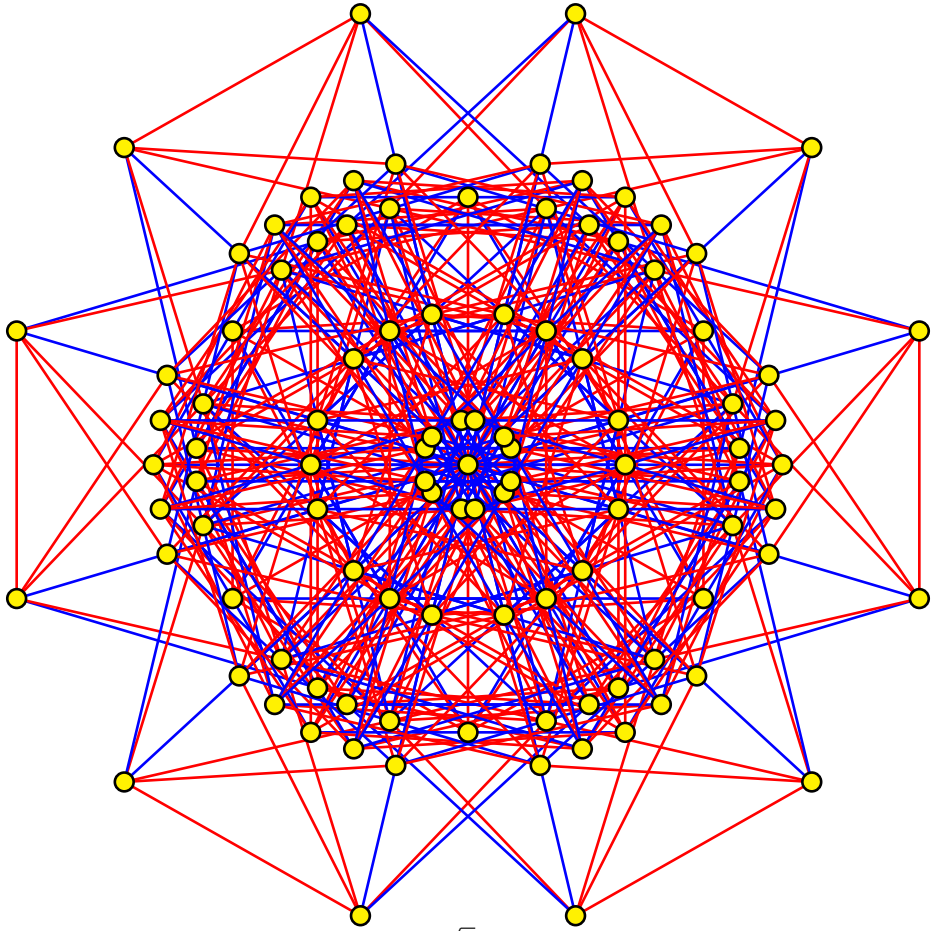}
\end{tikzpicture}
\caption{A 103-vertex $\{1, 2/\sqrt{3}\}$-graph with chromatic number 5.}
\label{g103fig}
\end{figure}

Of course, this graph is much too large to be handled directly. It has $312$ unit edges, $177$ edges of length $2/\sqrt{3}$, and chromatic number $5$. Despite its size, Sage takes only a couple of minutes to verify this. A list of vertices is available at the url \cite{website}.
\end{proof}

\section{\bf Conclusions and directions of future research}

Suppose one wants to color the plane such that no two points distance 1 or $d$ apart can be colored identically. We proved that at least five colors are needed
for $d$ taking any of the following values:
\begin{equation*}
\frac{\sqrt{5}+1}{2},\,\, \sqrt{3},\,\, \frac{\sqrt{6}+\sqrt{2}}{2},\,\, \frac{1}{2}\sqrt{3^{1/4}\cdot 2\sqrt{2}+2\sqrt{3}+2},\,\,\sqrt{3/2+\sqrt{33}/6},\,\, \frac{\sqrt{5}}{\sqrt{3}},\,\, 2, \,\,\frac{2}{\sqrt{3}}
\end{equation*}

In a couple of instances we proved results of the following form:
\begin{equation*}
\text{If}\,\, \chi(\mathbb{E}^2, \{1, d\})=4\,\,\text{then}\,\,  \chi(\mathbb{E}^2, \{1,d, d'\})=4.
\end{equation*}

This allowed us to reduce showing that $\chi(\mathbb{E}^2, \{1, d\})\ge 5$ to proving that $\chi(\mathbb{E}^2, \{1,d, d'\})\ge 5$ instead.

As already mentioned, if $d=\sqrt{11/3}$ one can prove an implication of the form:
\begin{equation*}
\text{If}\,\, \chi(\mathbb{E}^2, \{1\})=4\,\,\text{then}\,\,  \chi(\mathbb{E}^2, \{1,d\})=4.
\end{equation*}

It is tempting to propose the following
\begin{conj}
There exist values $d\neq 1$ such that
\begin{equation*}
\text{If}\,\, \chi(\mathbb{E}^2, \{1\})=5\,\,\text{then}\,\,  \chi(\mathbb{E}^2, \{1,d\})=5.
\end{equation*}
\end{conj}
If one can then find a $6$-chromatic $\{1,d\}$-graph the above conjecture would immediately imply that $\chi(\mathbb{E}^2)\ge 6$.
At this time however, no such graphs are known.

\section{\bf Appendix}
\begin{figure}[htp]
\centering
\includegraphics[width=\linewidth]{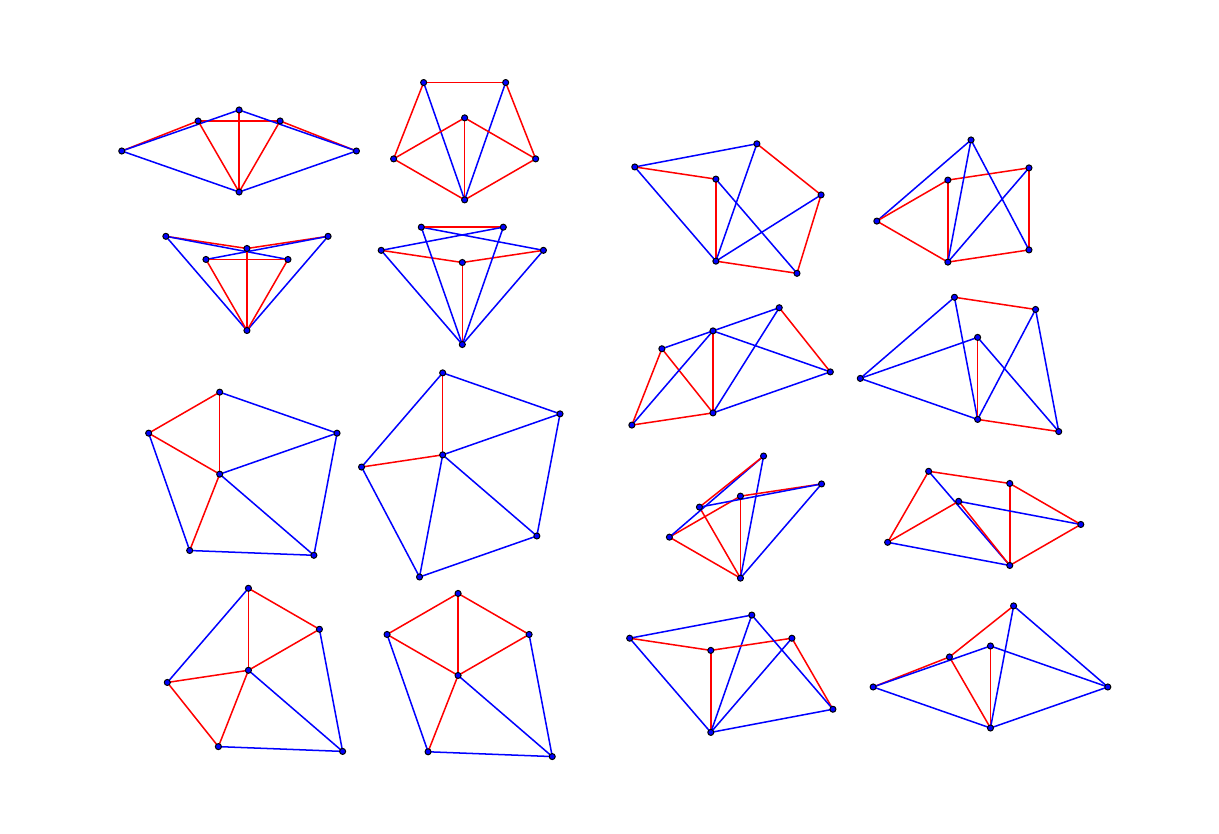}
\vspace{-0.5cm}
\caption{Sixteen different embeddings of the wheel graph $W_6$ as a $\{1, d\}$-graph. Here $d=\frac{1}{2}\sqrt{3^{1/4}\cdot 2\sqrt{2}+2\sqrt{3}+2}$.}
\label{16embeddings}
\end{figure}

\end{document}